\documentclass[12pt]{amsart}
 \usepackage{amsmath,amsfonts,amssymb}
 \usepackage{amsthm}
 \usepackage{array}
 \usepackage{graphics}
 \usepackage{amsmath}
 \usepackage{amssymb}
 \usepackage{verbatim}

 \usepackage{t1enc}
 \usepackage[utf8]{inputenc}

 \newtheorem{thm}{Theorem}
 \newtheorem{lem}{Lemma}
 \newtheorem{cor}{Corollary}

 \newtheorem*{rem}{Remark}

\newcommand{\kk}{\mathbb{K}}

\usepackage{amsmath}

\title{Some notes on Erd\H{o}s' last equation}

\author{Csaba~S\'andor}
\email{sandor.csaba@ttk.bme.hu}
\address{Department of Stochastics, Institute of Mathemetics, Budapest University of Technology and Economics, M\H{u}egyetem rkp. 3., H-1111 Budapest, Hungary; \newline \hspace*{4mm}
HUN-REN Alfr\'ed R\'enyi Institute of Mathematics, Re\'altanoda utca 13--15., H-1053 Budapest, Hungary
MTA--HUN-REN RI Lend\"ulet ``Momentum'' Arithmetic Combinatorics Research Group, Re\'altanoda utca 13--15., H-1053 Budapest, Hungary}

\author{Maciej Zakarczemny}
\email{mzakarczemny@pk.edu.pl}
\address{Department of Applied Mathematics, Cracow University of Technology, Warszawska 24, 31-155 Kraków, Poland.}

\thanks{
Csaba S\'andor is supported by the NKFIH Grants No. K129335, KKP 144059 and the Lend\"ulet ''Momentum'' program of the Hungarian Academy of Sciences (MTA)}

\begin{document}
\begin{abstract}
We study the Diophantine equation $$n^h(x_1 + \ldots + x_n) = x_1 x_2 \cdots x_n,$$ where $h \in \mathbb{Z}_{\ge 0}$, and $x_1 \le \ldots \le x_n \in \mathbb{Z}^+$. We define the solution set $ELE(h, n)$ consisting of all such integer tuples. We derive bounds on $x_n$ and $x_{n-1}$ in $ELE(h,n)$. We prove that
$x_1 \cdot\ldots\cdot x_n \le n^h(n^h + 1)(n^h + n - 1)$.
We also show that for $h = 1$,
$x_{n-1} \le n + \max_{d \mid 2n^2 - 2n,\, d \le \sqrt{2n^2 - 2n}} d$,
and for $h \ge 2$, the inequality $x_{n-1} \le 2n^h$ holds in $ELE(h, n)$.
Moreover, we study the minimal value $g(h, n)$ of the largest coordinate $x_n$, and define $g(h) = \min_{n \ge 2} g(h,n)$. We prove that
$$
\liminf_{n \to \infty} \frac{g(h,n)}{\frac{\log \log n \cdot \log \log \log \log n}{\log \log \log n}} \ge 1.
$$
We also estimate the growth of the number of solutions. We show that for any $\varepsilon > 0,$ we have $\sum_{m \le n} |ELE(h,m)| < n^{h+1+\varepsilon}
$
for sufficiently large $n$, and that $|ELE(h,n)| > n^{1 - \varepsilon}$
for infinitely many $n$.
\end{abstract}



\maketitle

\section{Introduction}

It is well known that the so-called "Equal Sum and Product Problem" concerns the number of solutions to the equation
\begin{equation}\label{ESP}
x_1+\ldots +x_n=x_1\cdot \ldots \cdot x_n,
\end{equation}
where $1 \le x_1\le x_2\le \dots \le x_n$, $x_i \in \mathbb{Z}$, $n\ge 2$.

On 19 August 1996, a month before his death, Paul Erd\H{o}s formulated the following problem in a letter written to Richard Guy (see \cite{Guy}):

"Do you know anything about the equation
\begin{equation}\label{ELE}
n(x_1+x_2+\ldots +x_n)=x_1x_2\cdot\ldots\cdot x_n,
\end{equation}
$x_1\le x_2\le \dots \le x_n$, where the $x_i$ are positive integers? Denote the number of solutions in $n$ by $f(n)$; $f(n)>n^{\varepsilon }$ for some $\varepsilon >0$ seems to be true? ..."

Equation (\ref{ELE}) is called "Erd\H{o}s' last equation".

Let $h$ be a non-negative integer. A generalization of the above equations is as follows:
\begin{equation}\label{ELEh}
n^h(x_1+x_2+\ldots +x_n)=x_1x_2\cdot\ldots\cdot x_n,
\end{equation}
where $x_1\le x_2\le \dots \le x_n$ and $x_i$ are positive integer, $n\ge 2$.\\
Let us define the set
\begin{equation}\label{ELE(h)}
\begin{aligned}
ELE(h,n) = \Big\{ (x_1, \dots ,x_n): \, &n^h (x_1 + x_2 + \ldots + x_n) = x_1 x_2 \cdot\ldots\cdot x_n, \\
&1 \leq x_1 \leq x_2 \leq \dots \leq x_n, x_i \in \mathbb{Z} \Big\}.
\end{aligned}
\end{equation}

In 1998, Kurlandchik and \cite{KN} proved

\begin{thm}\label{KN}
For equation (\ref{ESP}), we have
$$
x_1\cdot \ldots\cdot x_n\le 2n \qquad \mathrm{ and} \qquad x_n\le n.
$$
\end{thm}

In 2019, Shiu \cite{Shiu} proved
\begin{thm}
For the equation (\ref{ELE}),
we have
$$
x_n\le n(2n-1) \qquad \mathrm {and} \qquad
\limsup_{n\to \infty }\frac{x_{n-1}}{n}\le 1+\sqrt{2}.
$$
\end{thm}

In the paper \cite{BFT}, the following equation is studied:
\begin{equation}
x_1 x_2 \cdot \ldots\cdot x_n = h(n)(x_1 + x_2 + \ldots + x_n),
\end{equation}
where \( n \) is a fixed positive integer, \( h(n) \) is an integer-valued arithmetic function, and the variables \( 1 \leq x_1 \leq x_2 \leq \cdots \leq x_n \) take positive integer values. The authors of \cite{BFT} obtained an upper bound on the number of solutions to this equation:

\begin{equation}
N(n) \leq \frac{1}{2} h(n) n (h(n) n + 1).
\end{equation}

We prove the following generalization of the above results.

\begin{thm}\label{ELEHXN}
Let $h$ be a non-negative integer.
If $(x_1,\dots ,x_n)\in ELE(h,n)$, then
\begin{equation}
x_n\le n^h(n^h+n-1)
\end{equation}
and equality holds iff $(x_1,\dots ,x_n)=(1,\dots ,1,n^h+1,n^h(n^h+n-1))$.
\end{thm}
\begin{proof}[Proof of Theorem \ref{ELEHXN}]
We can assume that the solution of $ELE(h,n)$ has $n-l-1$ unit elements, where
$1 \le l \le n-1$.
Unit elements are $x_1 = x_2 = \ldots = x_{n-l-1} = 1$, while the non-unit elements are
$2 \leq x_{n-l} \leq x_{n-l+1} \leq \ldots \leq x_n.$
\begin{comment}
We have:
$
2^l x_n \leq x_{n-l} \cdot x_{n-l+1} \cdot \ldots \cdot x_n = n^h (x_{1} + x_{2} + \ldots + x_n) \leq n^{k+1} x_n.
$
Thus, $1 \leq l \leq (h+1) \log_2 n.$\\
\end{comment}
Let us denote
$$
P = \prod_{i=n-l}^{n-1} x_i\geq 2^l\ge 2, \quad S = n-1-l+\sum_{i=n-l}^{n-1} x_i\geq n+l-1 \ge n.
$$
Thus, we have $P x_n = \prod_{i=1}^n x_i,\,S + x_n = \sum_{i=1}^n x_i$ which from equation \eqref{ELE} gives:
$$n^h (S + x_n)=Px_n.$$
Put $w = P - n^h,$ hence $wx_n = n^h S.$
%Therefore $w|n^kS.$
%We also have $w = \tfrac{n^h S}{x_n},$ thus $1 \le w \le \tfrac{n^h S}{2}.$

Let us denote
$R = \frac{S-n+l+1}{P}=\frac{x_{n-l} + \ldots + x_{n-1}}{x_{n-l} \cdot \ldots \cdot x_{n-1}}.$
Since $l\ge 1$, we have
$$
R= \frac{x_{n-l} + \ldots + x_{n-1} }{x_{n-l} \cdot \ldots \cdot x_{n-1}} \le\frac{l x_{n-1}}{x_{n-l} \cdot \ldots\cdot x_{n-2} \cdot x_{n-1}} \le \frac{l}{2^{l-1}} \leq 1.
$$

We have
$
x_{n-l} + \ldots + x_{n-1} = R P = R (n^h + w)$ and
$
S = n - 1 - l + R (n^h + w).
$
Since $l \geq 1$, $w \ge 1$, $R \leq 1$ we get:
\begin{align*}
x_n = &\frac{n^h S}{w} = \frac{1}{w} n^h (n - 1 - l + R n^h) + R n^h \leq \\ &n^h (n - 2+ n^h) + n^h = n^h (n^h + n - 1).
\end{align*}
The equality holds iff $l=1,\,w=1,\,R=1,$ equivalent to $(x_1,\dots ,x_n)=(1,\dots ,1,n^h+1,n^h(n^h+n-1))$.
\end{proof}

\begin{thm}\label{ELEHPROD}
Let $h$ be a non-negative integer. If $(x_1,\dots ,x_n)\in ELE(h,n)$, then
\begin{equation}
x_1\cdot \ldots \cdot x_n\le n^h(n^h+1)(n^h+n-1)
\end{equation}
and equality holds iff $(x_1,\dots ,x_n)=(1,\dots ,1,n^h+1,n^h(n^h+n-1))$.
\end{thm}

\begin{proof}
By Theorem \ref{KN}, we may assume $h\ge 1$.
In the notation from the proof of Theorem \ref{ELEHXN}, we have:
$$x_1x_2\cdot\ldots\cdot x_n=Px_n=n^h(S+x_n)=n^h(S+\tfrac{n^hS}{w})=n^h(1+\tfrac{1}{w}n^h)S.$$
Hence
$$x_1x_2\cdot\ldots\cdot x_n=n^h(1+\tfrac{1}{w}n^h)(n-l-1+R(n^h+w)).$$
Since $l\ge 1,\,R\le 1$ we get
$$x_1x_2\cdot\ldots\cdot x_n\le n^h(1+\tfrac{1}{w}n^h)(n^h+n+w-2).$$
Note that
\begin{align*}
w=&\frac{n^hS}{x_n}=n^h\frac{n-1-l+\sum\limits_{i=n-l}^{n-1}x_i}{x_n}\le n^h\frac{n-1-l+lx_n}{x_n}\le \\ &n^h\frac{n-1-l}{x_n}+n^hl< n^{h}(n-1-l)+ n^hl=n^h(n-1).
\end{align*}
%But $l\le (h+1)\log_2n$, thus
%\begin{equation}\label{boundh}
%w\le\tfrac{1}{2}n^{h+1}-n^h+(h+1)n^h\log_2n.
%\end{equation}
The inequality
$$n^h(1+\tfrac{1}{w}n^h)(n^h+n+w-2)\le n^h(n^h+1)(n^h+n-1),$$
which is the equivalent of inequality
$$0\le (1-\tfrac{1}{w})(n^{2h}+n^{h+1}-2n^h-w),$$
follows from the bounds $1\le w<n^h(n-1)\le n^{2h}+n^{h+1}-2n^h$. The~equality holds iff $l=1,\,w=1,\,R=1.$
\end{proof}

\begin{thm}\label{thn1h} Let $n\ge 2$. The following bounds hold on
$x_{n-1}$:
\begin{enumerate}
\item For the equation (\ref{ELE}) (i.e. $h=1$) we have $$x_{n-1}\le n+\max_{d|2n^2-2n, d\le \sqrt{2n^2-2n}}d$$ and equality holds iff \begin{align*}&(x_1,\dots ,x_{n-2},x_{n-1},x_{n})=\\ &(1,\dots ,1,n+\max_{d|2n^2-2n, d\le \sqrt{2n^2-2n}}d,n+\min_{d|2n^2-2n, d\ge \sqrt{2n^2-2n}}d).\end{align*}
\item For equation (\ref{ELEh}) and $h\ge 2$, if $(x_1,\dots ,x_n)\in ELE(h,n)$, we have $x_{n-1}\le 2n^h$ and the equality holds iff $$(x_1,\dots ,x_n)=(1,\dots ,1,2n^h,2n^h+n-2).$$
\end{enumerate}
\end{thm}
\begin{proof}
Let $h$ be a positive integer. Assume that the solution of $ELE(h,n)$ has $n-1-l$ unit elements, where $1 \le l \le n-1$.\\
If $l=1$, then we have two non-unit elements. Hence
\begin{equation}\label{eqhn1}n^h(n^h+n-2)=(x_n-n^h)(x_{n-1}-n^h)\ge(x_{n-1}-n^h)^2.
\end{equation}
Therefore $$x_{n-1}\le n^h+\sqrt{n^h(n^h+n-2)}.$$
If $h=1$, then $x_{n-1}-n\le\sqrt{2n^2-2n}$ and $x_{n-1}-n|2n^2-2n$.
Thus
$$x_{n-1}-n\le \max_{d|2n^2-2n, d\le \sqrt{2n^2-2n}}d.$$
The equality holds iff $$(x_1,\dots ,x_{n-2},x_{n-1},x_{n})=$$ $$(1,\dots ,1,n+\max_{d|2n^2-2n, d\le \sqrt{2n^2-2n}}d,n+\min_{d|2n^2-2n, d\ge \sqrt{2n^2-2n}}d).$$\\
For $h\ge 2$ (the case $h=2, n=2$ is verified directly), it follows from \eqref{eqhn1}, the arithmetic–geometric mean inequality, the inequalities $$n^h(n^h+n-2)>(n^h+\tfrac{1}{2}n-\tfrac{3}{2})^2,\qquad x_{n-1}\le x_n,$$ and $x_{n-1},x_n\in\mathbb{Z}^+,$ that $(x_{n-1}-n^h)+(x_n-n^h)\ge 2n^h+n-2$. For $(x_{n-1}-n^h)+(x_n-n^h)= 2n^h+n-2$, we have $x_{n-1}=2n^h$ and $x_n=2n^h+n-2$. If $(x_{n-1}-n)+(x_n-n)> 2n^h+n-2$, then we have $x_{n-1}<2n^h$. Hence
$$x_{n-1}\le 2n^h$$
and the equality holds iff $(x_1,\dots ,x_n)=(1,\dots,1,2n^h,2n^h+n-2)$.\\
Assume that $l\ge 2$. Let us denote
$$
A = \prod_{i=n-l}^{n-2} x_i\geq 2^{l-1}\ge 2, \quad B = n-1-l+\sum_{i=n-l}^{n-2} x_i\geq n+l-3 \ge n-1.
$$
Thus, we have $A x_{n-1} x_n = \prod_{i=1}^n x_i,\,\,B +x_{n-1}+x_n = \sum_{i=1}^n x_i$ which from equation \eqref{ELEh} gives:
$$n^h (B + x_{n-1}+ x_n)=Ax_{n-1}x_n.$$
Let us denote
$T = \frac{B-n+1+l}{A}=\frac{x_{n-l} + \ldots + x_{n-2}}{x_{n-l} \cdot \ldots \cdot x_{n-2}}.$
Since $l\ge 2$, we have
\begin{equation*}
T= \frac{x_{n-l} + \ldots + x_{n-2} }{x_{n-l} \cdot \ldots \cdot x_{n-2}} \le\frac{(l-1) x_{n-2}}{x_{n-l} \cdot \ldots\cdot x_{n-3} \cdot x_{n-2}} \le \frac{l-1}{2^{l-2}} \le 1.\end{equation*}
Put $g = Ax_n - n^h,$ hence $ n^h (B+x_n)=gx_{n-1}.$\\
We have $x_{n-l} + \ldots + x_{n-2} = T A = \tfrac{T}{x_n}(n^h + g),$ hence
$$
B = n - 1 - l + \tfrac{T}{x_n}(n^h + g),$$
\begin{align*}
x_{n-1}=&\tfrac{n^h}{g}(B+x_n)=\tfrac{n^h}{g}(n - 1 - l + \tfrac{T}{x_n}(n^h + g)+x_n)=\\
&\tfrac{n^{2h}}{gx_n}T+n^h\tfrac{x_n}{g}+n^h\tfrac{1}{g}(n - 1 - l)+n^h\tfrac{T}{x_n}.
\end{align*}
Since $x_n=\tfrac{1}{A}(n^h+g)$ we get
\begin{equation}\label{eqhn2}
x_{n-1}=\tfrac{n^{2h}}{gx_n}T+n^h\tfrac{1}{A}\tfrac{n^h}{g}+
n^h\tfrac{1}{A}+n^h\tfrac{1}{g}(n - 1 - l)+n^h\tfrac{T}{x_n}.
\end{equation}
If $x_n< 2n^h,$ then $x_{n-1}\le x_n< 2n^h$. Thus, we may assume that $x_n\ge 2n^h$ and $g=Ax_n-n^h\ge 2x_n-n^h\ge  3n^h.$\\
Since $l \geq 2$, $A \ge 2$, $g \ge 3n^h$, $x_n\ge 2n^h$, $T \le 1$ we get:
%$$x_{n-1} <\tfrac{n^{2k}}{3n^k2n^k}2+n^k\tfrac{1}{2}\tfrac{n^k}{3n^k}+\tfrac{1}{2}n^k+n^k\tfrac{1}{3n^k}(n-3)+n^k2\tfrac{1}{2n^k}$$
$$x_{n-1}\le \tfrac{1}{6}+\tfrac{1}{6}n^h+\tfrac{1}{2}n^h+\tfrac{1}{3}n+\tfrac{1}{2}<2n^{h}.$$
\end{proof}
The previous theorem shows that the maximum of $x_{n-1}$ is taken for $(x_1,\dots ,x_n)\in ELE(h,n)$, $h\ge 1$, when $x_i=1$ for $1\le i\le n-2$. This is not true for $h=0$. Let $a$ be a~positive integer. Then $
(a^2+a-1,2a-1)=1.$ It follows that there are infinitely many primes $p$ of the form $p=(a^2+a-1)k-(2a-1)$, where $k$ is a~positive integer. Let $n=p+1=(a^2+a-1)k-2a+2$. Then for $(1,...,1,x_{n-1},x_n)\in ELE(0,n)$ we have $(x_{n-1}-1)(x_n-1)=n-1=p$, that is $x_{n-1}=2$.\\
On the other hand, for $a>1$, if we take $x_1=\dots =x_{n-3}=1$, $x_{n-2}=a$, $x_{n-1}=a+1$ and $x_n=k$, then we have $$x_1+\ldots +x_n=(a^2+a-1)k-2a-1+a+(a+1)+k=a(a+1)k=x_1\cdot \ldots\cdot x_n.$$
It seems to be a hard problem to determine the maximum of $x_{n-1}$ in $(x_1,\dots ,x_n)\in ELE(0,n)$. We simply prove the following weaker result.

\begin{thm} For equation (\ref{ESP}) we have $x_{n-1}\le 1+\sqrt{n-1}$ and equality holds if and only if $n-1$ is a square number and $$(x_1,\dots ,x_{n-2},x_{n-1},x_{n})=(1,\dots 1,1+\sqrt{n-1},1+\sqrt{n-1}).$$
\end{thm}

\begin{proof}
For $ n < 5 $, the only solutions are $ (2,2), (1,2,3), (1,1,2,4)$. So, the inequality holds. We can assume that $ n \geq 5 $.
The proof is analogous to the proof of Theorem \ref{thn1h}.
If $l=1$, then from equation \eqref{eqhn1} for $h=0$, we have
$$n-1=(x_n-1)(x_{n-1}-1)\ge (x_{n-1}-1)^2.$$
Hence, $x_{n-1}\le 1+\sqrt{n-1}$ and equality holds if and only if $n-1$ is a~square and $$(x_1,\dots ,x_{n-2},x_{n-1},x_{n})=(1,\dots 1,1+\sqrt{n-1},1+\sqrt{n-1}).$$
If $l\ge 2$, then similarly to equation \eqref{eqhn2}, for $ h=0 $, we have that
$$x_{n-1}=\tfrac{1}{gx_n}T+\tfrac{1}{Ag}+\tfrac{1}{A}+\tfrac{1}{g}(n-1-l)+\tfrac{T}{x_n},$$
where $g=Ax_n-1,\,g\in\mathbb{Z},\,l\ge 2, A\ge 2, T\le 1.$\\
If $x_n\le \sqrt{n-1}$, then $x_{n-1}\le x_n<1+\sqrt{n-1}.$\\
If $x_n>\sqrt{n-1},$ then $g>2\sqrt{n-1}-1$ and since $n\ge 5$, we get
\begin{align*}x_{n-1}<&\tfrac{2}{(2\sqrt{n-1}-1)\sqrt{n-1}}+\tfrac{1}{2(2\sqrt{n-1}-1)}+\tfrac{1}{2}+\tfrac{1}{2\sqrt{n-1}-1}(n-3)+\tfrac{1}{\sqrt{n-1}}<\\ &1+\sqrt{n-1}.\end{align*}
\end{proof}
\begin{rem}
For equation (\ref{ESP}), Theorem~6 provides the bound $x_{n-1}\le 1+\sqrt{n-1},$ with equality attained exactly when $n-1$ is a square.
For other values of $n$, the exact maximum of $x_{n-1}$ is not known.
\end{rem}

In the next theorem, we consider as small $x_n$ as possible. For given $h\ge 1$ and $n\ge 2$, denote by $g(h,n)$ the least possible value of $x_n$ for equation (\ref{ELEh}). We define the function \( g(h) \), which minimizes \( g(h, n) \) for \( n \geq 2 \), that is, $\displaystyle g(h)=\min _{n\ge 2}\{ g(h,n) \}$.

\begin{thm}\label{MINXN} Let us suppose that $h\ge 1$.
\begin{enumerate}
\item  Then, we have $g(h)\ge 3$.
\item Then we have $$\frac{\log x}{\log \log x}\ll |\{ h:h\le x, g(h)=3  \} |\ll \log ^3x,   $$ for every $x\ge 10$.
\end{enumerate}
\end{thm}

\begin{proof}

(1): We prove by contradiction. Let us suppose that for some $h\ge 1$ and $n\ge 2$ we have $(x_1,\dots ,x_n)\in ELE(h,n)$ with $x_n=2$. Then the product $x_1\cdot\ldots\cdot x_n$ must be a power of two. It follows that $n=2^{k}$ for some nonnegative integer $k$ and $x_1+\ldots +x_{2^k}$ is also a power of two with $2^k<x_1+\ldots +x_{2^k}\le 2^{k+1}$. Hence $x_1+\ldots +x_{2^k}=2^{k+1}$, that is $x_i=2$ for every $1\le i\le 2^k$. Thus the equation (\ref{ELEh}) is equivalent to $2^{kh}2^{k+1}=2^{2^k}$, that is $2^{k}\equiv 1 \pmod{k}$. It follows that $k$ must be an odd integer. For $k=1$, we have $h=0$. Let us suppose that $k\ge 3$. Let $p$ be the least prime divisor of $k$. By Fermat's little theorem we have $2^{p-1}\equiv 1 \pmod{p}$ and $2^{k}\equiv 1\pmod{p}$. Hence $2^{(p-1,k)}\equiv 1\pmod{p}$, that is $2\equiv 1\pmod{p}$, a contradiction.\\
(2): First we prove the lower bound. Let $p$ be an odd prime number. By Fermat's little theorem we have $2^{p-1}\equiv 1 \pmod{p}$, thus  $h_p=\frac{2^{p-1}-p-1}{p}$ is an integer. If $p>3$, then $h_p>0$. For $x_1=\ldots =x_{2^{p-1}+1}=1$, $x_{2^{p-1}+2}=\ldots =x_{2^p-1}=2$ and $x_{2^p}=3$, we have $n=2^p$, $x_1+\ldots +x_n=3\cdot 2^{p-1}$ and $x_1\cdot\ldots\cdot x_n=3\cdot 2^{2^{p-1}-2}$, so $(x_1,\ldots ,x_{2^p})\in ELE(h_p,2^p)$. It is easy to see that for $p\ne p'$ we have $h_p\ne h_{p'}$. It follows that $$|\{ h:h\le x, g(h)=3  \} |\ge \pi (\log _2x)-2,$$ where $\pi (.)$ denote the prime counting function. By the Prime Number Theorem we have
$$|\{ h:1\le h\le x, g(h)=3  \} |\gg \frac{\log x}{\log \log x}.$$
Finally we prove the upper bound. First we verify that if $h\le x$ and $(x_1,\ldots ,x_n)\in ELE(h,n)$ with $x_n=3$, then $n\le x^2$ if $h$ and $x$ are large enough. If $x_n=3$, then $n=2^{a}3^{b}$ for some nonnegative integers $a$ and $b$. Similarly, $n<x_1+\ldots +x_n=2^{a'}3^{b'}$ for some nonnegative integers $a'$ and $b'$.
Tijdeman \cite[Theorem 1, Corollary]{Tij} proved that $$x_1+\ldots +x_n-n>\frac{n}{(\log n)^C},$$ for some positive effectively computable constant $C$. It follows
that the $n$-tuple $(x_1,\ldots ,x_n)$ contains either at least $\frac{n}{4(\log n)^C}$ $x_i$s equal to $2$ or at least $\frac{n}{6(\log n)^C}$ $x_i$s equal to $3$. Therefore $$x_1\cdot\ldots\cdot x_n\ge \min \{ 2^{\frac{n}{4(\log n)^C}},3^{\frac{n}{6(\log n)^C}} \} =2^{\frac{n}{4(\log n)^C}}.$$ Hence
$$3n^{h+1}\ge n^h(x_1+\ldots +x_n)=x_1\cdot\ldots\cdot x_n\ge 2^{\frac{n}{4(\log n)^C}},$$ that is $h\ge c_1\frac{n}{(\log n)^{C+1}} $ for some $c_1>0$. It follows that for $(x_1,\dots ,x_n)\in ELE(h,n)$, $x_n=3$ we have $n\le c_2h(\log h)^{C+1}$ for some $c_2>0$. If $h$ is large enough then $c_2h(\log h)^{C+1}\le h^2\le x^2$, that is $n\le x^2$.\\
Since $n=2^{a}3^b\le x^2$, then $a\le 2\log _2x$ and $b\le 2\log _3x$, that is for $(x_1,\dots ,x_n)\in ELE(h,n)$ with $x_n=3$ and $h\le x$ we have at most $c_3(\log x)^2$ choices for the integer $n$. Let $x_1+\ldots +x_n=m$. Then $n<m\le 3n$ and $m=2^{a'}3^{b'}$. For a given $c\le 2\log _2x+\log_23$ there is a unique integer $d$ such that $n<2^c3^d\le 3n$, that is the number of possible $m$ is at most $c_4\log x$. Hence, it is enough to prove that for any pair $(n,m)$ we have at most one $n$-tuple $(x_1,\dots ,x_n)$, with $x_n=3$ such that $x_1+\ldots +x_n=m$ and $(x_1,\dots ,x_n)\in ELE(h,n)$ for some positive $h$. Let us suppose that we have different $n$-tuples $(x_1,\dots ,x_n)$ and $(x_1',\dots ,x_n')$ with $x_n=3$ and $x_n'=3$ such that $x_1+\ldots +x_n=x_1'+\ldots +x_n'=m$ and $(x_1,\dots ,x_n)\in ELE(h,n)$ and $(x_1',\dots ,x_n')\in ELE(h',n)$ for some positive integers $h$ and $h'$. Let us suppose that $h\le h'$. If $h=h'$, then $$x_1\cdot\ldots\cdot x_n=n^hm=n^{h'}m=x_1'\cdot\ldots\cdot x_n'$$ implies, that $(x_1,\dots ,x_n)=(x_1',\dots ,x_n')$, a contradiction. If $h<h'$, then $x_\cdot \ldots\cdot x_n=n^hm=2^c3^d$ and $x_1'\cdot \ldots\cdot x_n'=n^{h'}m=2^c3^dn^{h'-h}=2^{c'}3^{d'}$ imply that $c\le c'$, $d\le d'$ and $c+d<c'+d'$. Then $$m=x_1+\ldots +x_n=1(n-(c+d))+2c+3d=n+c+2d$$ and $$m=x_1'+\ldots +x_n'=1(n-(c'+d'))+2c'+3d'=n+c'+2d',$$ a contradiction.
\end{proof}
\noindent For $n>e^{e^e}$, we introduce the notation
$$
L(n)=\frac{\log\log n \cdot \log\log\log n}{\log\log\log\log n},
$$
in order to simplify the expressions.\\
Luca, Pink and S\'{a}ndor \cite{LPS} proved that $$\liminf_{n\to \infty}\frac{g(1,n)}{L(n)}\ge 1.$$ The next theorem is a generalization of the above result.

\begin{thm}
Fix the positive integer $h$. Consider equation (\ref{ELEh}). Then, we have $$\liminf_{n\to \infty}\frac{g(h,n)}{L(n)}\ge 1.$$
\end{thm}

Let $\eta$ be an algebraic number of degree $d$ whose
minimal polynomial over the integers is $$f(x) = a_0 \prod_{i=1}^d (x - \eta_i)$$with $a_0\ge 1$.
The {\it absolute logarithmic height} of $\eta$ is defined as
$$
h(\eta):= \frac{1}{d}\left( \log a_0 + \sum_{i=1}^d \log \max \{|\eta_i|,1\}\right).
$$

Let $\kk \subset \mathbb{R}$ be an algebraic number field and let $D$ be the degree of
the field $\kk$. Let $\eta_1, \eta_2, \ldots, \eta_l$ be nonzero elements of $\kk$ and $b_1, \ldots, b_l$ be integers.
We put
\begin{equation}\label{eq:Lambda}
\Lambda = \prod_{i=1}^l \eta_i^{b_i} -1,
\end{equation}
and
\begin{equation} \label{eq:B}
B  \ge \max\{|b_1|, \ldots, |b_l|\}.
\end{equation}
Let $A_1, \ldots, A_l$ be positive integers such that $$A_j \geq h'(\eta_j) := \max \{Dh(\eta_j), |\log \eta_j|, 0.16
\}\quad {\text{\rm for}}\quad j\in\{1,\ldots , l\}.
$$

The  following is contained in Theorem 9.4 of \cite{BMS} which is a consequence of a result due to Matveev \cite{Matv}.

\begin{lem} \label{lem:Matveev}
If $\Lambda \neq 0$ then
\begin{equation*}
\label{ineq:matveev} \log |\Lambda| > -C(l,D)(1+\log{B})A_1A_2\cdot \ldots\cdot A_l,
\end{equation*}
where
$$
C(l,D):=1.4 \cdot 30^{l+3}l^{4.5}D^2(1+ \log D).
$$
\end{lem}

\begin{proof}
Let $2=p_2<p_2<\dots <p_t$ be the first $t$ prime numbers. Denote by $S_t$ the set of all positive integers composed of $p_1,p_2,\dots ,p_t$ and arrange this in increasing order
$$
  1=n_1<n_2<\dots
$$
We are going to give a lower bound for $\frac{n_{i+1}}{n_i}-1 $. We may assume that $n_{i+1}\le 2n_i$. Write $n_i=p_1^{a_1}p_2^{a_2}\cdot \ldots\cdot p_t^{a_t}$ and $n_{i+1}=p_1^{b_1}p_2^{b_2}\cdot \ldots\cdot p_t^{b_t}$. Then
$$
\frac{n_{i+1}}{n_i}-1=p_1^{a_1-b_1}p_2^{a_2-b_2}\cdot \ldots\cdot p_t^{a_t-b_t}-1.
$$
Observe that
$$
max_{1\le k\le t}|b_k-a_k|\le \max\{ a_k,b_k\}\le 2\log n_{i+1}\le 4\log n_i.
$$
Apply Lemma \ref{lem:Matveev} with $l=t$, $D=1$, $A_1=A_2=\dots =A_t=\log p_t$ and $B=4\log n_i$ to conclude that
\begin{equation}\label{eq1}
\begin{split}
\frac{n_{i+1}}{n_i}-1\ge \frac{1}{\exp \{1.4\cdot 30^{t+4}t^{4.5}(1+\log (4\log n_i))(\log p_t)^t \} }.
\end{split}
\end{equation}

Assume that there exist infinitely many $(x_1,\dots ,x_n)\in ELE(h,n)$ such that
$$
x_n\le L(n).
$$
Let us suppose that
$$
1=x_1=x_2=\dots =x_{n-l},\quad 2\le x_{n-l+1}\le \dots \le x_n.
$$
Then $n^h\cdot n\cdot L(n) \ge 2^l$, that is
$$
l\le (h+1)\log_2n+\log_2\log\log n +\log_2 \log\log\log\log n \le (h+2)\log _2n
$$
for large enough $n$. Then
$$
n<x_1+x_2+\ldots +x_n\le n+(h+2)\log_2 n \cdot L(n).
$$
Let $p_t$ be the largest prime number dividing $x_1x_2\cdot \ldots\cdot x_n$. Then\\ $x_1x_2\cdot \ldots\cdot x_n\in S_t$. It follows that $n,x_1+x_2+\ldots + x_n\in S_t$.  Clearly,
$$
\frac{x_1+x_2+\ldots +x_n}{n}-1\le (h+2)\frac{\log _2n}{n}L(n)
$$
By (\ref{eq1}) we have
\begin{align*}
&\frac{1}{\exp \{1.4\cdot 30^{t+4}t^{4.5}(1+\log (4\log n))(\log p_t)^t \} } \le \\ &(h+2)\frac{ L(n) \cdot  \log_2n}{n}\le \frac{1}{\sqrt{n}}
\end{align*}
if $n$ is large enough, that is
$$
0.5 \log n < 1.4\cdot 30^{t+4}t^{4.5}(1+\log (4\log n))(\log p_t)^t.
$$
By the Prime Number Theorem we have $\log p_t =(1+o(1))\log t$ as $t\to \infty$.
%$$
%\log n <\exp\{(1+o(1))t\log\log t\}
%$$
%$$
%\log \log n< (1+o(1))t\log\log t.
%$$
A simple calculation shows that $t>(1+o(1))\frac{\log\log n}{\log\log\log\log n}$.\\
By the Prime Number Theorem, we have $p_t>(1+o(1))L(n)$.\\ Since $p_t\mid x_1x_2\cdot \ldots\cdot x_n$ and $x_1\le x_2\le \dots \le x_n$ we have $$x_n\ge p_t\ge (1+o(1))L(n),$$ which completes the proof.

\end{proof}

%\begin{thm} Let us suppose that $h\ge 1$ and $n\ge 2$.
%\begin{enumerate}
%\item  For any $(x_1,\dots, x_n)\in ELE(h,n)$ we have $x_n\ge 3$.
%\item There exist infinitely many tuples $(h,n)$ and $(x_1,\dots, x_n)\in ELE(h,n)$ with $x_n=3$.
%\item If we fix a~positive integer $h$, then as $n\to \infty$, we have $%%$x_n\ge (1+o(1))\frac{\log \log n \cdot \log \log \log n}{\log \log \log %\log n}\ge 1.$$
%\end{enumerate}
%\end{thm}

%\noindent--------------------------------------\\
%NUMERICAL RESULTS (Istv\'{a}n Pink)\\
%--------------------------------------\\\\

In 2012, Weingartner \cite{AW} proved the following theorem.
\begin{thm}\label{Weingartner}
We have
$$
\sum_{m=2}^n|ELE(0,m)|=(1+o(1))\frac{ne^{2\sqrt{\log n}}}{2\sqrt{\pi}(\log n)^{3/4}}.
$$
\end{thm}

It is easy to see that for distinct $(x_1,\dots ,x_n),(x_1',\dots ,x_n')\in ELE(h,n)$, where $x_1 = x_2 = \ldots = x_{n-l-1} = 1$, $2 \leq x_{n-l} \leq x_{n-l+1} \leq \ldots \leq x_n$ and $x_1' = x_2' = \ldots = x_{n-l'-1}' = 1$, $2 \leq x_{n-l}' \leq x_{n-l+1}' \leq \ldots \leq x_n'$, we have $(x_{n-l},\dots ,x_n)\ne (x_{n-l'}',\dots ,x_n')$. Let $a_n$ denote the number of factorizations of the integer $n$ into factors larger than $1$, where the order of factors does not count. By Theorem \ref{ELEHPROD} we get $|ELE(h,n)|\le \sum_{m\le n^h(n^h+1)(n^h+n-1)}a_m$. In 1927, A. Oppenheim proved that $$\sum_{n\le x}a_n=(1+o(1))\frac{xe^{2\sqrt{\log x}}}{2\sqrt{\pi} (\log x)^{3/4}}$$ as $x\to \infty$. Let $h\ge 1$ and $\varepsilon >0$. It follows that $$\sum_{m\le n}|ELE(h,m)|\le n^{3h+\varepsilon}$$ if $n$ is large enough.

The next theorem improves this bound for $h\ge 1$.
%\begin{thm}\label{f(n)}
%\item Let $h\ge 1$. Then we have $|ELE(h,n)|<n^{2h+2}$ if $n$ is large enough.
%\end{thm}
%\begin{thm}\label{ELE(1,n)}
%\item Let $\varepsilon >0$. Then we have $\displaystyle \sum_{m\le n}|%ELE(1,m)|<n^{2+\varepsilon }$ if $n$ is large enough.
%\end{thm}
%\begin{proof}
%Let $m,n\in \mathbb{Z}^{+}$, $m\le n$. Let us suppose that for %\begin{equation*}m(x_1+\dots +x_m)=x_1\dots x_m.\end{equation*}
%Without loss of generality, we may assume $x_1=\dots =x_{m-l-1}=1$ and  $2\le x_{m-l}\le \dots \le x_m$. Then $x_1\dots x_{m-1}x_{m}\le m\cdot m %x_m$. Hence $x_{m-l}\dots x_{m-1}\le m^2\le n^2$.  We will show that for %given $(y_1,y_2,\ldots,y_l)$ such that  $2\le y_1\le \dots \le y_l$, %%%%%$y_i\in \mathbb{Z}^+$, $y_1\dots y_l\le n^2$ we have at most $n^{o(1)}$ %tuples $(y,m)$, $y,m\in \mathbb{Z}^+$, $m\le n$ with
%\begin{equation*}
%m(m-l-1+y_1+\dots +y_l+y)=y_1\dots y_ly.
%\end{equation*} Here \( o(1) \) refers to the limit as \( n \to \infty %\).
%Let $U=y_1+\dots +y_l-l-1$ and $V=y_1\dots y_l$, thus $U\le V\le n^{2}$. %Then
%$$y=\frac{m^2+mU}{V-m}=-U-V-m+\frac{(U+V)V}{V-m}.$$
%It follows that $V-m\mid (U+V)V$. Note that  $(U+V)V=O(n^4)$, , then the number of divisors satisfies $d((U+V)V)=n^{o(1)}$. Hence
$$
%\sum_{m\le n}|ELE(1,n)|\le n^{o(1)}|\{ (y_1,\dots y_l): y_i\ge 2, y_i\in \mathbb{Z}^+, y_1\dots y_l\le n^2 \} |.
%$$
%Oppenheim proved that
$$
%|\{ (y_1,\dots y_l): y_i\ge 2, y_i\in \mathbb{Z}^+, y_1\dots  y_l\le n^2 %\} |=n^{2+o(1)},
%$$
%which completes the proof.
%\end{proof}
\begin{thm}\label{ELE(h,n)}
Let $h\in\mathbb{Z}^+$ and $\varepsilon >0$.\\ Then we have $\displaystyle \sum_{m\le n}|ELE(h,m)|<n^{h+1+\varepsilon }$ if $n$ is large enough.
\end{thm}

\begin{proof}
Let $m,n\in \mathbb{Z}^{+}$, $m\le n$. Let us suppose that for \begin{equation*}m^h(x_1+\ldots +x_m)=x_1\cdot \ldots\cdot x_m.\end{equation*}
Without loss of generality, we may assume $x_1=\dots =x_{m-l-1}=1$ and $2\le x_{m-l}\le \dots \le x_m$. Then $x_1\cdot\ldots\cdot x_{m-1}x_{m}\le m^h\cdot m x_m$. Hence $x_{m-l}\cdot\ldots\cdot x_{m-1}\le m^{h+1}\le n^{h+1}$.
We will show that for given $(y_1,y_2,\ldots,y_l)$ such that $2\le y_1\le \dots \le y_l$, $y_i\in \mathbb{Z}^+$, $y_1\cdot\ldots\cdot y_l\le n^{h+1}$ we have at most $n^{o(1)}$ tuples $(y,m)$, $y,m\in \mathbb{Z}^+$, $m\le n$ with
\begin{equation*}
m^h(m-l-1+y_1+\ldots +y_l+y)=y_1\cdot\ldots\cdot y_ly.
\end{equation*} Here \( o(1) \) refers to the limit as \( n \to \infty \).

Let $U=y_1+\ldots +y_l-l-1$ and $V=y_1\cdot\ldots\cdot y_l$. %Note that $(y_1-1)(y_2-1)\dots(y_l-1)\ge 1$, thus by induction
Clearly, $U\le V\le n^{h+1}$.
Then
\begin{equation*}y=\frac{m^{h+1}+m^hU}{V-m^h}=-m-U+\frac{UV+mV}{V-m^h}.\end{equation*}
It follows that $$mV \equiv -UV \pmod{V-m^h}$$ and $$0\equiv (UV+mV)^h\equiv \sum_{i=1}^{h}\binom{h}{i}(UV)^i(mV)^{h-i}+(mV)^h \equiv$$
$$
\sum_{i=1}^{h}\binom{h}{i}(UV)^i(-UV)^{h-i}+V^{h+1}\equiv V^{h}(V+U^h\sum_{i=1}^{h}\binom{h}{i}(-1)^{h-i})\equiv
$$
$$
V^h(V-(-1)^hU^h) \pmod{V-m^h}
$$
%$V-m^h\mid UV+mV$, thus  $V-m^h\mid (UV+mV)^h$.
%Hence, $V-m^h\mid V^{h+1}+(UV)^h\sum\limits_{i=0}^{h-1}(-1)^i$

Note that  $|V^h(V-(-1)^hU^h)|\le 2n^{2(h+1)^2}$.
Let $d(n)$ denote the number of positive divisors of an integer $n$. It is known that for any $\varepsilon>0,$ we have
$
d(n)=O(n^{\varepsilon}),\quad\text{as } n\to\infty.
$
Hence $d(n)=n^{o(1)}$ as $n\to\infty$. Let us suppose that $V\ne (-1)^hU^h$. Since $V^h(V-(-1)^hU^h)=O(n^{2(h+1)^2})$ and $h$ is fixed (independent of $n$) , the number of divisors satisfies:
$$
d(V^h(V-(-1)^hU^h))=n^{2(h+1)^2o(1)}=n^{o(1)}.
$$
If $V=(-1)^h U^h$, then $h=2h'$ is an even positive integer and $V=U^h$. We have
$$
y=-m+U+\frac{U^{2h'}(U+m)}{U^h-m^h}=$$ $$-m+U+\frac{U^{2h'}(U+m)}{(U^2-m^2)(U^{2h'-2}+U^{2h'-4}m^2+\ldots +m^{2h'-2})}=
$$
$$
-m+U+\frac{U^{2h'}}{U^{2h'-1}-U^{2h'-2}m+U^{2h'-3}m^2+-\ldots -m^{2h'-1}}.
$$
It follows that $$U^{2h'-1}-U^{2h'-2}m+U^{2h'-3}m^2+-\ldots -m^{2h'-1}\mid U^{2h'}. $$
So the number of possible $m$ is at most $(2h'-1)d(U^{2h'})=n^{o(1)}$.
Hence
$$
\sum_{m\le n}|ELE(h,n)|\le n^{o(1)}|\{ (y_1,\dots y_l): y_i\ge 2, y_i\in \mathbb{Z}^+, y_1\cdot\ldots\cdot y_l\le n^{h+1} \} |.
$$
Oppenheim \cite{Opp} proved that
$$
|\{ (y_1,\dots y_l): y_i\ge 2, y_i\in \mathbb{Z}^+, y_1\cdot\ldots\cdot  y_l\le n^{h+1} \} |=n^{h+1+o(1)},
$$
which completes the proof.
\end{proof}
A direct consequence of Theorem \ref{Weingartner} is the following.
\begin{cor}
Let $\varepsilon >0$. Then $|ELE(0,n)|<n^{1+\varepsilon}$ if $n$ is large enough.
\end{cor}
The next theorem shows that this bound is sharp.
\begin{thm}\label{lower}
Let $h$ be a non-negative integer, and let $\varepsilon >0$. Then $|ELE(h,n)|>n^{1-\varepsilon}$ for infinitely many positive integers $n$.
\end{thm}
\begin{proof}%[Proof of Theorem \ref{lower}]
It is enough to prove that there exist infinitely many $n$ such that there are at least $n^{1-\varepsilon}$ $n$-tuples $(x_1,\dots ,x_n)$, where there exists an $i$ with $x_i=x_{i+1}=\dots =x_{i+h-1}=n$ and $x_1+\ldots +x_{i-1}+nh+x_{i+h}+\ldots +x_n=x_1\cdot \ldots \cdot x_{i-1}x_{i+h}\cdot \ldots \cdot x_n$.
Let $m$ be a~positive integer. Denote by $p_m$ the $m$-th prime number. Let $p_u$ be the largest prime less than or equal to $h+1$. If
$h=0$, we take $u=0$. Let $P_m^{(h)}=p_{u+1}p_{u+2}\cdot \ldots\cdot p_{u+m}$,\, $P_0^{(h)}=1$. By the Prime Number Theorem,
%$p_m=(1+o(1))m\log m$ as $m\to \infty$. Hence,
\begin{equation}\label{e6}
P_m^{(h)}=e^{\sum_{k=u+1}^{u+m}\log p_k}=e^{\sum_{k=2}^m(\log k +\log\log k +o(1))}.
\end{equation}
By the Euler–Maclaurin summation formula (see \cite{dB,Nat}),
%We have (see also \cite[Theorem 6.3]{Nat}), \begin{equation}\label{e7} \begin{split}\sum_{k=2}^m\log k=&\sum_{k=2}^m\left(\int_k^{k+1} \log xdx +O\left(\log(k+1)-\log(k)\right)\right) \\=&\int_2^{m+1}\log xdx+O(\log m)\\=&m\log m-m+O(\log m).\end{split}\end{equation}Analogously, \begin{equation}\label{e8}\begin{split}\sum_{k=2}^m\log\log k=&\sum_{k=2}^m\left(\int_k^{k+1} \log\log xdx +O\left(\log\log(k+1)-\log\log k\right)\right) \\=&\int_2^{m+1}\log\log xdx+O(\log\log m)\\=&\int_2^{m+1}\log\log x+\frac{1}{\log x}dx-\int_2^{m+1}\frac{1}{\log x}dx+O(\log\log m)\\=&m\log\log m+O\left(\frac{m}{\log m}\right).\end{split}\end{equation}
%By \eqref{e6},\,\eqref{e7},\,\eqref{e8} we obtain
\begin{equation}\label{pm}
P_m^{(h)}=e^{m\log m +m\log \log m -m +o(m)}\,\,\mathrm{as}\,\,m\to \infty.
\end{equation}
One may notice that the following equality holds:
\begin{equation}
\frac{P_{(i+1)m}^{(h)}}{P_{im}^{(h)}}=e^{m\log m+o(m\log m)}
\end{equation}
uniformly for $i=O(\log m)$ and $m\to \infty$.\\
Let us suppose that $k\ge 3$. By (\ref{pm})
$$
P_{k\lfloor \log k \rfloor }^{(h)}=e^{k\log^2 k +(2+o(1))k\log k\log \log k}.
$$
as $k\to\infty $. Let
\begin{equation}M_k=\max_{0\le i\le \lfloor \log k \rfloor }\frac{P_{(i+1)k}^{(h)}}{P_{ik}^{(h)}}=e^{k\log k +o(k\log k)}.\end{equation}
Denote by $B_n$ be the number of partitions of a set of size $n$, which are called Bell numbers. We will show that
\begin{equation}\label{lowerbound}
\sum_{l=1}^{h+1}\sum_{lP_{k\lfloor \log k \rfloor}^{(h)} -\lfloor \log k\rfloor M_k<n\le lP_{k\lfloor \log k \rfloor}^{(h)}}|ELE(h,n)|\ge B_k^{\lfloor \log k\rfloor }.
\end{equation}
It is known (see \cite{dB}) that $B_n=e^{n\log n-n\log\log n -n +o(n)}$ as $n\to \infty$. Hence
\begin{gather*}
\max_{\substack{1\le l\le h+1\\lP_{k\lfloor \log k \rfloor}^{(h)} -\lfloor \log k\rfloor M_k<n\le lP_{k\lfloor \log k \rfloor}^{(h)}}}|ELE(h,n)|\\ \ge
\frac{B_k^{\lfloor \log k\rfloor }}{(h+1)\lfloor \log k \rfloor M_k} = e^{k\log^2 k-(1+o(1))k\log k\log\log k}.
\end{gather*}
Let $n_k$ be the integer, where the maximum is taken. Note that $$n_k\le (h+1)P_{k\lfloor\log k\rfloor}^{(h)}<P_{(k+1)\lfloor \log (k+1) \rfloor}^{(h)} -\lfloor \log (k+1)\rfloor M_{k+1}\le n_{k+1},$$ if $k$ is sufficiently large. Thus the set $\{n_k:k\ge 2\}$ is infinite.
We have also
$$\log n_k=(1+o(1))\log P_{k\lfloor \log k\rfloor }^{(h)}=(1+o(1))k\log^2 k,$$
$$
\log\log n_k=(1+o(1))\log k,
$$
and
$$
\log\log\log n_k=(1+o(1))\log\log k.
$$
Hence
\begin{align*}
|ELE(h,n_k)|&\ge e^{k\log^2k-(1+o(1)) k\log k\log\log k }\\
&=e^{(1-o(1))k\log ^2k}\\%&\ge n_ke^{-(3-o(1))\frac{\log n_k\log\log\log n_k}{\log\log n_k}}\\
&=n_k^{1-o(1)}.
\end{align*}
Finally we prove (\ref{lowerbound}). For a~positive integer $k$, let $$\{ S_1^{(1)},\ldots ,S_{b_1}^{(1)} \},\,\{ S_1^{(2)},\ldots ,S_{b_2}^{(2)} \},\ldots,\,\{ S_{1}^{(B_k)},\ldots ,S_{b_{B_k}}^{(B_k)} \}$$ be the partitions of $\{ 1,...,k\}.$\\
For a given $(i_1,\ldots ,i_{\lfloor \log k\rfloor })\in \{ 1,\ldots ,B_k\} ^{\lfloor \log k\rfloor}$ and any $1\le j\le \lfloor \log k\rfloor $, define
$$
y_{b_{i_1}+\ldots +b_{i_{j-1}}+v}=\prod_{s\in S_v^{(i_j)}}p_{(j-1)k+s},
$$
where $1\le v\le b_{i_j}$. For any $1\le j\le \lfloor \log k\rfloor $,
$$
\prod_{1\le v\le b_{i_j}}y_{b_{i_1}+\ldots +b_{i_{j-1} }+v}=\frac{P_{jk}^{(h)}}{P_{(j-1)k}^{(h)}}
$$
and
$$
\sum_{1\le v\le b_{i_1}+\ldots + b_{i_{\lfloor \log k \rfloor }}}y_v\le M_k\lfloor \log k\rfloor .
$$
We complete the sequence  $(y_1,\ldots ,y_{b_{1}+\ldots +b_{\lfloor \log k\rfloor }})$ with an $h+1$, with an integer $1\le l\le h+1$ such that $$h+1\mid y_1+y_2+\ldots +y_{b_{1}+\ldots +b_{\lfloor \log k\rfloor }}+l-(b_{1}+\ldots +b_{\lfloor \log k\rfloor }-1),$$ with $h$ integer $$\frac{(h+1)lP_{k\lfloor \log k \rfloor } ^{(h)}-(y_1+y_2+\ldots +y_{b_{1}+\ldots +b_{\lfloor \log k\rfloor }}+l-(b_{1}+\ldots +b_{\lfloor \log k\rfloor })-1)}{h+1}$$ and with
\begin{align*}\frac{(h+1)lP_{k\lfloor \log k \rfloor } ^{(h)}-(y_1+y_2+\ldots +y_{b_{1}+\ldots +b_{\lfloor \log k\rfloor }}+l-(b_{1}+\ldots +b_{\lfloor \log k\rfloor })-1)}{h+1}-\\ (b_{1}+\ldots +b_{\lfloor \log k\rfloor })-2-h\end{align*} 1s to get the sequence $(x_1,\dots ,x_n)$ with $$n=\frac{(h+1)lP_{k\lfloor \log k \rfloor }^{(h)}-(y_1+y_2+\ldots +y_{b_{1}+\ldots +b_{\lfloor \log k\rfloor }}+l-(b_{1}+\ldots +b_{\lfloor \log k\rfloor })-1)}{h+1}.$$ Then $n^h(x_1+\ldots +x_n)=x_1x_2\cdot\ldots\cdot x_n $ with
$$
lP_{k\lfloor \log k \rfloor }^{(h)}-M_k\lfloor \log k\rfloor \le n\le lP_{k\lfloor \log k \rfloor }^{(h)}
$$

Note that, if $$(i_1,\ldots ,i_{\lfloor \log k\rfloor })\ne (i_1',\ldots ,i_{\lfloor \log k\rfloor }'),$$ then
$$
(y_1,\ldots ,y_{b_{i_1}+\ldots b_{i_{\lfloor \log k\rfloor }}})\ne (y_1',\ldots ,y_{b_{i_1'}+\ldots b_{i_{\lfloor \log k\rfloor }'}}'),
$$
so different sequences $(i_1,\ldots ,i_{\lfloor \log k\rfloor })\in \{ 1,\ldots ,B_k\} ^{\lfloor \log k\rfloor}$ yield different solutions to equation \eqref{ELEh}, which proves the lower bound \eqref{lowerbound}.
\end{proof}

\end{document}